\documentclass{amsart}
\newtheorem{thm}{Theorem}
\newtheorem{cor}[thm]{Corollary}
\newtheorem{lem}[thm]{Lemma}

\newtheorem{defn}[thm]{Definition}

\numberwithin{equation}{section}

\begin{document}

\title[Smooth Approximation of Lipschitz functions]{Smooth Approximation
of Lipschitz functions on Riemannian manifolds}
\author{D. Azagra, J. Ferrera, F. L\'{o}pez-Mesas, Y. Rangel}

\address{Departamento de An\'{a}lisis Matem\'{a}tico\\ Facultad de
Matem\'{a}ticas\\ Universidad Complutense\\ 28040 Madrid, Spain}

\date{January 31, 2006}

\email{azagra@mat.ucm.es, ferrera@mat.ucm.es,
FLopez\_Mesas@mat.ucm.es}




\begin{abstract}
We show that for every Lipschitz function $f$ defined on a
separable Riemannian manifold $M$ (possibly of infinite
dimension), for every continuous $\varepsilon:M\to (0,+\infty)$,
and for every positive number $r>0$, there exists a $C^\infty$
smooth Lipschitz function $g:M\to\mathbb{R}$ such that
$|f(p)-g(p)|\leq\varepsilon(p)$ for every $p\in M$ and
$\textrm{Lip}(g)\leq\textrm{Lip}(f)+r$. Consequently, every
separable Riemannian manifold is uniformly bumpable. We also
present some applications of this result, such as a general
version for separable Riemannian manifolds of
Deville-Godefroy-Zizler's smooth variational principle.
\end{abstract}

\maketitle

\section{Introduction and main results}

It is well known, and very useful, that every Lipschitz function
$f:\mathbb{R}^{d}\to\mathbb{R}$ can be uniformly approximated by
$C^\infty$ smooth Lipschitz functions whose Lipschitz constants
are the same as $f$'s. This can be done easily by considering the
integral convolutions
$$f_{n}=\int_{\mathbb{R}^{d}}f(y)\varphi_{n}(x-y)dy,$$
where the $\varphi_{n}$ satisfy
$\int_{\mathbb{R}^{d}}\varphi_{n}=1$ and
$\textrm{supp}(\varphi_{n})\subset B(0, 1/n)$. This method of
smooth approximation has many advantages over other standard
procedures like smooth partitions of unity, as the integral
convolutions preserve many geometrical properties that $f$ may
have, such as convexity or Lipschitzness. Indeed, if $f$ is
$K$-Lipschitz then $f_{n}$ is $K$-Lipschitz as well.

For finite-dimensional Riemannian manifolds, Greene and Wu
\cite{Greene5, Greene3, Greene4} used a refinement of this
integral convolution procedure to get very useful results on
smooth approximation of convex or Lipschitz functions defined on
Riemannian manifolds (in fact they applied this method to prove
several theorems about the structure of complete noncompact
manifolds of positive curvature). It turned out, however, that,
when one is interested in approximating a convex function $h$ by
$C^\infty$ convex functions this method works out (that is, gives
convex $f_{n}$'s) in Riemannian manifolds only when the function
$h$ is {\em strictly convex} (see \cite{Greene3}); and also that,
when one needs to perform a $C^{0}$-fine approximation of a
Lipschitz function $f$ by $C^\infty$ smooth Lipschitz functions,
the approximations $f_{n}$ have Lipschitz constants which are
arbitrarily close to the Lipschitz constant of $f$ (but are not
equal in general).

Unfortunately, the integral convolution method breaks down in
infinite dimensions (due to the lack of a suitable measure like
Lebesgue's one), and other methods have to be employed instead. It
is well known that $C^\infty$ smooth partitions of unity exist on
every Riemannian manifold and can of course be used to get
$C^0$-fine approximation of continuous functions by $C^\infty$
smooth functions. On the other hand, Moulis \cite{Moulis} showed
that $C^1$-fine approximations of $C^1$ smooth functions by
$C^\infty$ smooth functions are also available on
infinite-dimensional Riemannian manifolds. We should also mention
that infimal convolutions with squared geodesic distances can be
used to regularize convex functions on Riemannian manifolds of
nonpositive sectional curvature \cite{AFconvinf}.

However, no one seems to have considered the natural question
whether every Lipschitz function $f$ defined on an
infinite-dimensional Riemannian manifold can be $C^0$-finely
approximated by $C^\infty$ smooth functions $g$ whose Lipschitz
constants also approximate the Lipschitz constant of $f$. We think
this is a very interesting question because many functions arising
from geometrical problems on Riemannian manifolds are Lipschitz
but not $C^1$ smooth (the distance function to a closed subset of
a manifold is a typical instance), so smooth approximations which
almost preserve Lipschitz constants can be very helpful in the
analysis of such problems.

In fact we were motivated to study this question by an open
problem from \cite{AFL2}: whether a version for Riemannian
manifolds of the Deville-Godefroy-Zizler smooth variational
principle \cite{DGZ2, DGZ} holds for every complete separable
Riemannian manifold. This is a very interesting problem because
the DGZ variational principle is an invaluable tool in the
(nonsmooth) analysis of Hamilton-Jacobi equations defined on
Riemannian manifolds. In \cite{AFL2} we were able to prove such a
variational principle under the assumption that the manifold was
{\em uniformly bumpable} (see Definition \ref{definition of
uniform bumpability} below), but the question whether or not every
Riemannian manifold is uniformly bumpable remained open.

In this note, as a consequence of our result on smooth Lipschitz
approximation we will answer these two questions in the
affirmative: every separable Riemannian manifold $M$ is uniformly
bumpable and, consequently, if $M$ is complete, satisfies the DGZ
smooth variational principle.

On the other hand, we have been informed that Garrido, Jaramillo
and Rangel \cite{GaJaRan} have recently established an
infinite-dimensional version of the Myers-Nakai theorem
\cite{Myers, Nakai} under the assumption that the manifold is
uniformly bumpable (and therefore their result holds in fact for
every separable Riemannian manifold). This encourages us to expect
that the result we present on smooth Lipschitz approximation (as
well as the fact that every separable Riemannian manifold is
uniformly bumpable) will find more applications beyond the DGZ
variational principle or the infinite-dimensional Myers-Nakai
theorem.

Let us now state the main result of this note.

\begin{thm}\label{smooth approximation of Lipschitz functions}
Let $M$ be a separable Riemannian manifold, let $f:M\to\mathbb{R}$
be a Lipschitz function, let $\varepsilon:M\to (0,+\infty)$ be a
continuous function, and $r>0$ a positive number. Then there
exists a $C^\infty$ smooth Lipschitz function $g:M\to\mathbb{R}$
such that $|f(p)-g(p)|\leq\varepsilon(p)$ for every $p\in M$, and
$\textrm{Lip}(g)\leq\textrm{Lip}(f)+r$.
\end{thm}
Here $\textrm{Lip}(f)$ and $\textrm{Lip}(g)$ stand for {\em the}
Lipschitz constants of $f$ and $g$, respectively, that is,
    $$
    \textrm{Lip}(f)=\inf\{ L\geq 0 : |f(p)-f(q)|\leq L d(p,q)\}=
    \sup\{\frac{|f(p)-f(q)|}{d(p,q)}: p, q\in M, p\neq q\}.
    $$
(recall that $f$ is said to be Lipschitz provided
$\textrm{Lip}(f)$ is finite, and if $L\geq \textrm{Lip}(f)$ then
we say that $L$ is {\em a} Lipschitz constant of $f$, or that $f$
is $L$-Lipschitz).

We should stress that we do not know whether a similar statement
holds for infinite-dimensional separable Banach spaces with
$C^\infty$ smooth bump functions, and that even in the case when
$M$ is a separable infinite-dimensional Hilbert space this result
seems to be new.

Theorem \ref{smooth approximation of Lipschitz functions} will be
proved in the next section. Let us now deduce the announced
consequences.

We first recall the definition of uniform bumpability given in
\cite{AFL2}.
\begin{defn}\label{definition of uniform bumpability}
A Riemannian manifold $M$ is {\em uniformly bumpable} provided
there exist numbers $R>1$ (possibly large) and $r>0$ (possibly
small) such that for every $p\in M$, $\delta\in (0, r)$ there
exists a $C^1$ smooth function $b:M\to\mathbb [0,1]$ such that:
\begin{enumerate}
\item $b(p)=1$ \item $b(q)=0$ if $d(q,p)\geq\delta$ \item
$\sup_{q\in M}\|db(q)\|_{q}\leq R/\delta$.
\end{enumerate}
\end{defn}
We have:
\begin{cor}\label{all manifolds are uniformly bumpable}
All separable Riemannian manifolds are uniformly bumpable. In fact
the constant $R$ in Definition \ref{definition of uniform
bumpability} can always be chosen to be any number bigger than 1,
the number $r$ any positive number, and the function $b$ of class
$C^\infty$.
\end{cor}
\begin{proof}
Let $R>1$, $0<\delta<r$, and $p\in M$ be given, and consider the
function $f:M\to [0,1]$ defined by
$$
f(q)=\left\{%
\begin{array}{ll}
    1-\frac{1}{\delta} d(q,p), & \textrm{ if }\, d(q,p)\leq\delta ; \\
    0, & \textrm{ if }\, d(q,p)\geq\delta. \\
\end{array}%
\right.
$$
It is clear that $f$ is $\frac{1}{\delta}$-Lipschitz and satisfies
$f(p)=1$, and $f=0$ off $B(p,\delta)$. By Theorem \ref{smooth
approximation of Lipschitz functions}, for any $\varepsilon>0$
there exists a $C^\infty$ smooth function $g:M\to\mathbb{R}$ such
that $|g(q)-f(q)|\leq\varepsilon$ for all $q\in M$ and
$\textrm{Lip}(g)\leq\frac{1}{\delta}+\varepsilon$. Now take a
$C^\infty$ smooth function $\theta:\mathbb{R}\to [0,1]$ such that
\begin{enumerate}
\item[{(i)}] $\theta(t)=0$ for $t\leq\varepsilon$; \item[{(ii)}]
$\theta(t)=1$ for $t\geq 1-\varepsilon$, and \item[{(iii)}]
$\textrm{Lip}(\theta)\leq\frac{1+\varepsilon}{1-2\varepsilon}$,
\end{enumerate}
and define $b(q)=\theta(g(q))$ for all $q\in M$. Then it is clear
that $b(p)=1$, $b(q)=0$ if $d(q,p)\geq\delta$, and
$$
\sup_{q\in M}\|db(q)\|_{q}=\textrm{Lip}(b)\leq
\textrm{Lip}(\theta)\textrm{Lip}(g)\leq
\frac{1+\varepsilon}{1-2\varepsilon}\left(\frac{1}{\delta}+\varepsilon\right)\leq
\frac{R}{\delta}
$$
if $\varepsilon$ is chosen small enough.
\end{proof}
As a consequence, the version of the Deville-Godefroy-Zizler
variational principle which was proved in \cite{AFL2} for
uniformly bumpable complete Riemannian manifolds is now seen to
hold for every complete separable Riemannian manifold.
\begin{cor}[DGZ smooth variational principle for Riemannian manifolds]\label{general
version of DGZ variational principle} Let $M$ be a complete
Riemannian manifold modelled on a separable Hilbert space, and let
$F:M\longrightarrow (-\infty,+\infty]$ be a lower semicontinuous
function which is bounded below, $F\not\equiv+\infty$. Then, for
each $\delta>0$ there exists a bounded $C^1$ smooth and Lipschitz
function $\varphi:M\longrightarrow\mathbb{R}$ such that:
\begin{enumerate}
\item $F-\varphi$ attains its strong minimum on $M$, \item
$\|\varphi\|_{\infty}:=\sup_{p\in M}|\varphi(p)|<\delta$, and
$\|d\varphi\|_{\infty}:=\sup_{p\in M}\|d\varphi(p)\|_{p}<\delta$.
\end{enumerate}
\end{cor}

\medskip

\section{Proof of the main Theorem}

The proof combines all of the most important approximation methods
we know of, that is: integral convolutions, partitions of unity,
and infimal convolutions. We first obtain, by using exponential
charts and infimal convolutions, local $C^1$ smooth Lipschitz
approximations of $f$, next we regularize these local
approximations by resorting to a result of Moulis's \cite{Moulis}
(which partially relies on the use of integral convolutions on
finite dimensional subspaces of the separable Hilbert space $X$ on
which $M$ is modelled), making sure that the Lipschitz estimates
are preserved, and finally we glue all the local approximations
together with the help of a specially constructed partition of
unity.

We begin with the precise statement of Moulis's result.
\begin{thm}[Moulis]\label{Moulis theorem}
Let $G$ be an open subset of a separable Hilbert space $X$, let
$f:G\to\mathbb{R}$ be a $C^1$ smooth function, and
$\varepsilon:G\to (0, +\infty)$ be a continuous function. Then
there exists a $C^\infty$ function $g:G\to \mathbb{R}$ such that
$|f(x)-g(x)|\leq\varepsilon(x)$ and
$\|f'(x)-g'(x)\|\leq\varepsilon(x)$ for every $x\in G$.
\end{thm}
For a proof see \cite{Moulis}, or, for a version of this theorem
that holds in every Banach space with an unconditional basis and a
smooth bump function, see \cite{AFGJL}.

\medskip

The first step in the proof of Theorem \ref{smooth approximation
of Lipschitz functions} is to prove a weaker statement in the
special case when $M=X$, for a constant $\varepsilon$, and
assuming $f$ is bounded. This can be done by combining Moulis's
theorem with Lasry and Lions's regularization technique of sup-inf
convolutions \cite{LasryLions}.
\begin{thm}\label{special case}
Let $(X, \|\cdot\|)$ be a separable Hilbert space, let
$f:X\to\mathbb{R}$ be a bounded and Lipschitz function, and let
$\varepsilon>0$. Then there exists a $C^\infty$ smooth Lipschitz
function $g:X\to\mathbb{R}$ such that $|f(x)-g(x)|\leq\varepsilon$
for every $x\in X$, and
$\textrm{Lip}(g)\leq\textrm{Lip}(f)+\varepsilon$.
\end{thm}
\begin{proof}
Let us denote $K=\textrm{Lip}(f)$. Because $f$ is Lipschitz and
bounded on $X$, according to the main Theorem of
\cite{LasryLions}, the functions
$$
x\mapsto (f_{\lambda})^{\mu}(x):= \sup_{z\in X}\inf_{y\in
X}\{f(y)+ \frac{1}{2\lambda}\|z-y\|^{2}-
\frac{1}{2\mu}\|x-z\|^{2}\}
$$
are of class $C^{1,1}$ on $X$ and converge to $f$ uniformly on $X$
as $0<\mu<\lambda\to 0$. So let us pick $\lambda$ and $\mu$ with
$0<\mu<\lambda$ and small enough so that
$$
|(f_{\lambda})^{\mu}(x)-f(x)|\leq\frac{\varepsilon}{2} \eqno(1)
$$
for all $x\in X$. We first see that $(f_{\lambda})^{\mu}$ is
$K$-Lipschitz on $X$. This is an immediate consequence of the fact
that the operations of inf- and sup- convolutions (with squared
norms or with any other kernel) preserve the Lipschitz constants
of the functions to be regularized: that is, if $h:X\to\mathbb{R}$
is $L$-Lipschitz on $X$ then the function
$$
h_{\lambda}(x)=\inf_{y\in X}\{h(y)+\frac{1}{2\lambda}\|x-y\|^{2}\}
$$
is $L$-Lipschitz on $X$ as well. Indeed, note first that
$$
\inf_{y\in X}\{h(y)+\frac{1}{2\lambda}\|x-y\|^{2}\}= \inf_{y\in
X}\{h(x-y)+\frac{1}{2\lambda}\|y\|^{2}\},
$$
so the function $h_{\lambda}$ can be redefined as
$$
h_{\lambda}(x')=\inf_{y\in
X}\{h(x'-y)+\frac{1}{2\lambda}\|y\|^{2}\}. \eqno(2)
$$
Now take two points $x, x'\in X$ and a number $\delta>0$. By the
definition of $\inf$ in $(2)$ we can find a point $y'\in X$ such
that
$$
h(x'-y')+\frac{1}{2\lambda}\|y'\|^{2}-\delta\leq\inf_{y\in
X}\{h(x'-y)+\frac{1}{2\lambda}\|y\|^{2}\},
$$
hence
\begin{eqnarray*}
& & h_{\lambda}(x)-h_{\lambda}(x')=\inf_{y\in
X}\{h(x-y)+\frac{1}{2\lambda}\|y\|^{2}\}-\inf_{y\in
X}\{h(x'-y)+\frac{1}{2\lambda}\|y\|^{2}\}\leq\\
& &
h(x-y')+\frac{1}{2\lambda}\|y'\|^{2}-h(x'-y')-\frac{1}{2\lambda}\|y'\|^{2}+\delta=\\
& & h(x-y')-h(x'-y')+\delta\leq L\|x-y'-(x'-y')\|+\delta=
L\|x-x'\|+\delta,
\end{eqnarray*}
that is $h_{\lambda}(x)-h_{\lambda}(x')\leq L\|x-x'\|+\delta$, and
by symmetry and by letting $\delta$ go to $0$ we get that
$|h_{\lambda}(x)-h_{\lambda}(x')|\leq L\|x-x'\|$, that is
$h_{\lambda}$ is $L$-Lipschitz. An identical proof shows that the
same is true of the function $h^{\mu}$ defined by
$$
h^{\mu}(x)=\sup_{z\in X}\{h(z)-\frac{1}{2\mu}\|x-z\|^{2}\}.
$$
Therefore the function $(f_{\lambda})^{\mu}$ has the same
Lipschitz constant as $f$, namely $K$.

Now, since $ (f_{\lambda})^{\mu}$ is $C^{1}$ smooth, we can use
Moulis's theorem to find a $C^\infty$ smooth function
$g:X\to\mathbb{R}$ such that
$$
|g(x)-(f_{\lambda})^{\mu}(x)|\leq\frac{\varepsilon}{2}, \textrm{
and } \|g'(x)-((f_{\lambda})^{\mu})'(x)\|\leq\varepsilon \eqno (3)
$$
for all $x\in X$. By combining $(1)$ and $(3)$ we obtain that
$|f(x)-g(x)|\leq\varepsilon$ and also
$$
\textrm{Lip}(g)=\sup_{x\in X}\|g'(x)\|\leq\sup_{x\in X}
\|((f_{\lambda})^{\mu})'(x)\|+\varepsilon\leq
K+\varepsilon=\textrm{Lip}(f)+\varepsilon.
$$
\end{proof}

In the proof of Theorem \ref{smooth approximation of Lipschitz
functions} we will have to use the fact that a locally
$K$-Lipschitz function defined on a Riemannian  manifold is is
globally $K$-Lipschitz.
\begin{lem}\label{locally K-Lipshitz is
equivalent to K-Lipschitz} Let $M$ be a Riemannian manifold and
let $f:M\to\mathbb{R}$ be a function which is locally
$K$-Lipschitz (that is, for every $a\in M$ there exists
$\delta=\delta(a)>0$ such that $|f(p)-f(q)|\leq K d(p,q)$ for all
$p,q\in B(a,\delta)$). Then $f$ is $K$-Lipschitz on $M$.
\end{lem}
\begin{proof}
This fact is proved, for instance, in \cite[Lemma 2]{Greene3} in
the setting of finite-dimensional Riemannian manifolds, but it is
clear that the same argument is also valid in the
infinite-dimensional case.
\end{proof}

\medskip

Let us start with the proof of Theorem \ref{smooth approximation
of Lipschitz functions}. In this proof $X$ will stand for the
separable Hilbert space on which the manifold $M$ is modelled, and
$B(p,\delta)$ will denote the open ball of center $p$ and radius
$\delta$ in $M$, that is $B(p,\delta)=\{q\in M: d(q,p)<\delta\}$.
We will also put $K=\textrm{Lip}(f)$ for short.

With no loss of generality, we can assume that $\varepsilon(p)\leq
r/2$ for all $p\in M$ (if necessary just replace $\varepsilon$
with the continuous function $p\mapsto \min\{\varepsilon(p),
r/2\}$). Also, let us fix any number $\varepsilon'>0$ small enough
so that
$$
(K(1+\varepsilon')+\varepsilon')(1+\varepsilon')<K+\frac{r}{2}.
$$
Now, for every $p\in M$, let us choose $\delta_{p}>0$ small enough
so that the exponential mapping is a bi-Lipschitz $C^\infty$
diffeomorphism of constant $1+\varepsilon'$ from the ball
$B(0_{p}, 3\delta_{p})\subset TM_{p}$ onto the ball $B(p,
3\delta_{p})\subset M$ (see \cite[Theorem 2.3]{AFL2}). Moreover,
by continuity of $f$ and $\varepsilon$, we can assume that the
$\delta_{p}$ also are sufficiently small so that
$\varepsilon(q)\geq \varepsilon(p)/2$ and
$|f(q)-f(p)|\leq\varepsilon(p)/2$ for every $q\in B(p,
3\delta_{p})$

Since $M$ is separable we can take a sequence $(p_n)$ of points in
$M$ such that
$$
M=\bigcup_{n=1}^{\infty}B(p_{n}, \delta_{n}),
$$
where we denote $\delta_{n}=\delta_{p_{n}}$, and also
$\varepsilon_{n}=\varepsilon(p_{n})$. Now, for each
$n\in\mathbb{N}$ define a function $f_{n}: B(0_{p_{n}},
3\delta_{n})\to\mathbb{R}$ by
$$
f_{n}(x)=f(\exp_{p_{n}}(x)),
$$
which is $K(1+\varepsilon')$-Lipschitz. We can extend $f_{n}$ to
all of $TM_{p_{n}}$ by defining
$$
\hat{f}_{n}(x)=\inf_{y\in B(0_{p_{n}},
3\delta_{p_{n}})}\{f_{n}(y)+K(1+\varepsilon')\|x-y\|_{p}\}
$$
It is well known and very easy to show that $\hat{f}_{n}$ is a
Lipschitz extension of $f_{n}$ to all of $TM_{p_{n}}$, with the
same Lipschitz constant $K(1+\varepsilon')$. The function
$\hat{f}_{n}$ is bounded on bounded sets (because it is Lipschitz)
but is not bounded on all of $TM_{p_{n}}$. Nevertheless we can
modify $\hat{f}_{n}$ outside the ball $B(0_{p_{n}}, 4\delta_{n})$
so as to make it bounded on all of $TM_{p_{n}}$. For instance, put
$C=\sup\{ |\hat{f}_{n}(x)|+1 : x\in B(0_{p_{n}}, 4\delta_{n})\}$,
and define $\tilde{f}_{n}:TM_{p_{n}}\to\mathbb{R}$ by
    $$
    \tilde{f}_{n}(x)=
  \begin{cases}
    -C & \text{ if } \hat{f}_{n}(x)\leq -C, \\
    \hat{f}_{n}(x) & \text{ if } -C\leq \hat{f}_{n}(x) \leq C, \\
    +C & \text{ if } C\leq \hat{f}_{n}(x).
  \end{cases}
    $$
It is clear that $\tilde{f}_{n}$ is bounded on all of $TM_{p_{n}}$
and has the same Lipschitz constant as $\hat{f}_{n}$, which is
less than or equal to $K(1+\varepsilon')$. That is,
$\tilde{f}_{n}$ is a bounded $K(1+\varepsilon')$-Lipschitz
extension of $f_{n}$ to $TM_{p_{n}}$.

\medskip

Next we are going to construct a $C^\infty$ smooth partition of
unity subordinated to the covering $\{B(p_{n},
2\delta_{n})\}_{n\in\mathbb{N}}$ of $M$ and to estimate the
Lipschitz constant of each of the functions of this partition of
unity. Let us take a $C^{\infty}$ smooth function
$\theta_{n}:\mathbb{R}\to [0,1]$ such that $\theta_{n}=1$ on
$(-\infty, \delta_{n}]$ and $\theta_{n}=0$ on $[2\delta_{n},
+\infty)$, and define
$$
\varphi_{n}(p)=\left\{%
\begin{array}{ll}
    \theta_{n}(\|\exp^{-1}_{p_{n}}(p)\|_{p_{n}}), & \hbox{
    if $p\in B(p_{n}, 3\delta_{n})$;} \\
    0, & \hbox{ otherwise.} \\
\end{array}%
\right.
$$
It is clear that each of the functions
$\varphi_{n}:M\to\mathbb{R}$ is $C^\infty$ smooth and Lipschitz,
and satisfies $\varphi_{n}=1$ on the ball $B(p_{n}, \delta_{n})$,
and $\varphi_{n}=0$ on $M\setminus B(p_{n}, 2\delta_{n})$.

Let us define the functions $\psi_{k}=:M\to [0,1]$ by
$$
\psi_{k}=\varphi_{k}\prod_{j<k}(1-\varphi_{j}).
$$
It is clear that $\psi_k$ is $C_k$-Lipschitz, where
$$C_{k}:=\sum_{j\leq k}\textrm{Lip}(\varphi_{j}),$$ and it
is easy to see that
\begin{enumerate}
\item For each $p\in M$, if $k=k(p)=\min\{j: p\in B(p_{j},
\delta_{j})\}$ then, because $1-\psi_{k}=0$ on $B(p_{k},
\delta_{k})$, we have that $B(p_{k}, \delta_{k})$ is an open
neighborhood of $p$ that meets only finitely many of the supports
of the functions $\psi_{\ell}$. Indeed,
$\textrm{supp}(\psi_{\ell})\cap B(p_{k}, \delta_{k})=\emptyset$
for all $\ell>k$, and $\textrm{supp}(\psi_{k})\subset B(p_{k},
2\delta_{k})$; \item $\sum_{k}\psi_{k}=1$;
\end{enumerate}
that is, $\{\psi_n\}_{n\in\mathbb{N}}$ is a $C^{\infty}$ smooth
partition of unity subordinated to the covering $\{B(p_{n},
2\delta_{n})\}_{n\in\mathbb{N}}$ of $M$.

\medskip

Now, according to Theorem \ref{special case} we can find a
$C^\infty$ smooth function $g_{n}:TM_{p_{n}}\to\mathbb{R}$ such
that
$$
|g_{n}(x)-\tilde{f}_{n}(x)|\leq
\frac{\varepsilon_{n}}{2^{n+2}\left(C_{n}+1\right)}, \eqno(4)
$$
for all $x\in TM_{p_{n}}$, and
$$
\textrm{Lip}(g_{n})\leq
\textrm{Lip}(\tilde{f}_{n})+\varepsilon'\leq
K(1+\varepsilon')+\varepsilon'. \eqno(5)
$$

We are ready to define our approximation $g:M\to\mathbb{R}$ by
$$
g(p)=\sum_{n}\psi_{n}(p)g_{n}(\exp_{p_{n}}^{-1}(p))
$$
for any $p\in M$. Observe that if $p\in B(p_{n}, 3\delta_{n})$,
because $\exp_{p_{n}}$ is a $C^\infty$ diffeomorphism from
$B(0_{p_{n}}, 3\delta_{n})$ onto $B(p_{n}, 3\delta_{n})$, the
expression $\psi_{n}(p)g_{n}(\exp_{p_{n}}^{-1}(p))$ is well
defined and is $C^\infty$ smooth on $B(p_{n}, 3\delta_{n})$. On
the other hand, if $p\notin B(p_{n},
2\delta_{n})\supset\textrm{supp}(\psi_{n})$ then $\psi_{n}(p)=0$.
So we will agree that, for any $p\notin B(p_{n}, 3\delta_{n})$,
the expressions $\psi_{n}(p)g_{n}(\exp_{p_{n}}^{-1}(p))$ and
$g_{n}(\exp_{p_{n}}^{-1}(p))$ both mean zero (whether or not
$\exp_{p_{n}}^{-1}(p)$ makes sense in this case). With these
conventions, since the $\psi_{n}$ form a $C^\infty$ smooth
partition of unity it follows that $g$ is well defined and is
$C^\infty$ smooth on $M$.

Let us see that $g$ and $\textrm{Lip}(g)$ approximate $f$ and
$\textrm{Lip}(f)$, respectively, as required.

Fix any $p\in M$, and let $k=k(p)$ be as in $(1)$ above, so that
we have $\psi_{\ell}=0$ on $B(p_{k}, \delta_{k})$ for all
$\ell>k$, and let us estimate $|f-g|$. To simplify the notation
let us denote $x_{m}=\exp_{p_{m}}^{-1}(p)\in TM_{p_m}$ (and
observe that this expression may well make no sense for many
$m$'s, but in such cases the corresponding $g_{m}(x_{m})$ have
been defined to be zero; in the following estimation we will also
understand that $\tilde{f}_{m}(x_{m})$ means zero if $p$ is
outside the ball $B(p_{n}, 3\delta_{n})$). We have
\begin{eqnarray*}
& & |g(p) -f(p)|= \left|\sum_{m\leq k}\psi_{m}(p)
g_{m}(\exp_{p_{m}}^{-1}(p)) -f(p)\right|= \\
& & \left| \sum_{m\leq k}\psi_{m}(p)\left[ g_{m}(x_{m})- f(p)
\right] \right|=  \left| \sum_{m\leq k}\psi_{m}(p)\left[
g_{m}(x_{m})-\tilde{f}_{m}(x_{m})\right]\right|\leq \\ & &
\sum_{m\leq k}\psi_{m}(p)
\frac{\varepsilon_{m}}{2^{m+2}\left(C_{m}+1\right)}\leq
\sum_{m\leq k}\psi_{m}(p) \frac{\varepsilon_{m}}{2}\leq
\sum_{m\leq k}\psi_{m}(p)\varepsilon(p)=\\ & &
\sum_{m}\psi_{m}(p)\varepsilon(p)=\varepsilon(p).
\end{eqnarray*}
Finally, let us check that $\textrm{Lip}(g)\leq K+r$. Since $g$ is
defined on a Riemannian manifold, according to Lemma \ref{locally
K-Lipshitz is equivalent to K-Lipschitz}, it is enough to show
that $g$ is locally $\left(K+r\right)$-Lipschitz. Take a point
$a\in M$, and define $k=k(a)=\min\{j: a\in B(p_{j},
\delta_{j})\}$, so that $\textrm{supp}(\psi_{\ell})\cap B(p_{k},
\delta_{k})=\emptyset$ for all $\ell>k$. Let also
$$\delta_{a}=\min\{\delta_{1}, ..., \delta_{k}, \delta_{k}-d(a,
p_{k})\}, $$ and
$$
F_{p,q}=\{m\in\{1, ..., k\}:  B(p_{m},
2\delta_{m})\cap\{p,q\}\neq\emptyset\}.$$

\medskip

We have that, if $p, q\in B(a, \delta_{a})$, then:
\begin{enumerate}
\item[{(i)}] For every $m\in\{1, ..., k\}$, we have that $p\in
B(p_{m}, 3\delta_{m})$ whenever $q\in B(p_{m}, 2\delta_{m})$; and
symmetrically $q\in B(p_{m}, 3\delta_{m})$ whenever $p\in B(p_{m},
2\delta_{m})$. Consequently, for every $m\in F_{p,q}$ we have that
$p, q\in B(p_{m}, 3\delta_{m})$; in particular, if $m\in F_{p,q}$,
then $x_{m}:=\exp_{p_{m}}^{-1}(p)$ and
$y_{m}:=\exp_{p_{m}}^{-1}(q)$ are well defined, and (by using
$(5)$ above), we also have
\begin{eqnarray*}
& & |g_{m}(x_{m})-g_{m}(y_{m})|\leq
(K(1+\varepsilon')+\varepsilon')
\|\exp_{p_{m}}^{-1}(p)-\exp_{p_{m}}^{-1}(q)\|_{p_{m}}\leq\\
& & (K(1+\varepsilon')+\varepsilon')(1+\varepsilon')d(p,q)
\end{eqnarray*}

$$ $$ (recall
that $\exp_{p_{m}}^{-1}:B(p_{m}, 3\delta_{m})\to B(0_{p_{m}},
3\delta_{m})$ is $(1+\varepsilon')$-Lipschitz). \item[{(ii)}] If
$m\in\mathbb{N}\setminus F_{p,q}$ then $\psi_{m}(p)=0=\psi_{m}(q)$
(because $\textrm{supp}(\psi_{m})\subset B(p_{m}, 2\delta_{m})$
and $\textrm{supp}(\psi_{\ell})\cap B(p_{k},
\delta_{k})=\emptyset$ for all $\ell>k$).
\end{enumerate}
Hence we have that, for $p,q\in B(a, \delta_{a})$ (with the
notation $x_{m}=\exp_{p_{m}}^{-1}(p)$ and
$y_{m}=\exp_{p_{m}}^{-1}(q)$),
\begin{eqnarray*}
& & g(p)=\sum_{m\in F_{p,q}}g_{m}(x_{m})\psi_{m}(p), \,\,\,\,\,
g(q)=\sum_{m\in F_{p,q}}g_{m}(y_{m})\psi_{m}(q), \\
& & 1=\sum_{m\in F_{p,q}}\psi_{m}(p)\, = \, \sum_{m\in
F_{p,q}}\psi_{m}(q), \textrm{ and }\\
& & |g_{m}(x_{m})-g_{m}(y_{m})|\leq
(K(1+\varepsilon')+\varepsilon')(1+\varepsilon')d(p,q) \, \textrm{
whenever } m\in F_{p,q}.
\end{eqnarray*}
Bearing in mind all these facts, considering $(4)$ above, and
using that $\psi_{m}$ is $C_{m}$-Lipschitz, we can estimate, for
every $p,q\in B(a, \delta_{a})$,
\begin{eqnarray*}
& & | g(p)-g(q)|=\left|\sum_{m\in F_{p,q}}g_{m}(x_{m})\psi_{m}(p)
-\sum_{m\in F_{p,q}}g_{m}(y_{m})\psi_{m}(q)\right|=\\
& & | \left( \sum_{m\in F_{p,q}}g_{m}(x_{m})\psi_{m}(p)
\right)-f(p)+f(p)-
\left(\sum_{m\in F_{p,q}}g_{m}(x_{m})\psi_{m}(q)\right)+ \\
& & + \left(\sum_{m\in F_{p,q}}g_{m}(x_{m})\psi_{m}(q)\right)
-\left(\sum_{m\in F_{p,q}}g_{m}(y_{m})\psi_{m}(q) \right) |=\\ & &
| \left( \sum_{m\in F_{p,q}}g_{m}(x_{m})\psi_{m}(p)
\right)-\left(\sum_{m\in
F_{p,q}}f(p)\psi_{m}(p)\right)+\left(\sum_{m\in
F_{p,q}}f(p)\psi_{m}(q)\right)
\\
& &-\left(\sum_{m\in F_{p,q}}g_{m}(x_{m})\psi_{m}(q)\right) +
\left(\sum_{m\in F_{p,q}}g_{m}(x_{m})\psi_{m}(q)\right)
-\left(\sum_{m\in F_{p,q}}g_{m}(y_{m})\psi_{m}(q) \right) |=\\
& &\left| \sum_{m\in F_{p,q}}\left[
\left(g_{m}(x_{m})-f(p)\right)\left(\psi_{m}(p)-\psi_{m}(q)\right)+\left(
g_{m}(x_{m})-g_{m}(y_{m}) \right)\psi_{m}(q) \right] \right|\leq \\
 & &\left| \sum_{m\in F_{p,q}}
\left(g_{m}(x_{m})-f(p)\right)\left(\psi_{m}(p)-\psi_{m}(q)\right)\right|+
\left|\sum_{m\in F_{p,q}}\left( g_{m}(x_{m})-g_{m}(y_{m})
\right)\psi_{m}(q)
\right| \leq \\
& & \sum_{m\leq
k}\frac{\varepsilon_{m}}{2^{m+2}\left(C_{m}+1\right)}C_{m}d(p,q)+
\sum_{m\leq
k}(K(1+\varepsilon')+\varepsilon')(1+\varepsilon')d(p,q)\psi_{m}(q)\leq \\
& & \sum_{m\leq k}\frac{\varepsilon(a)}{2^{m+1}}d(p,q)+
\sum_{m\leq
k}(K(1+\varepsilon')+\varepsilon')(1+\varepsilon')d(p,q)\psi_{m}(q)\leq \\
& & \varepsilon(a)d(p,q)+(K(1+\varepsilon')+\varepsilon')(1+\varepsilon')d(p,q)\leq\\
& & \frac{r}{2}d(p,q)+\left(K+\frac{r}{2}\right)d(p,q)=
\left(K+r\right) d(p,q).
\end{eqnarray*}
This shows that $g$ is locally
$\left(\textrm{Lip}(f)+r\right)$-Lipschitz and concludes the proof
of Theorem \ref{smooth approximation of Lipschitz functions}. \qed

\medskip



\begin{thebibliography}{}

\bibitem{AFconvinf}
D. Azagra and J. Ferrera, {\em Inf-convolution and regularization
of convex functions on Riemannian manifolds of nonpositive
curvature}, preprint, 2005.

\bibitem{AFL2} D. Azagra, J. Ferrera, F. L\'{o}pez-Mesas, {\em
Nonsmooth analysis and Hamilton-Jacobi equations on Riemannian
manifolds}, J. Funct. Anal. 220 (2005) no. 2, 304-361.

\bibitem{AFGJL}
D. Azagra, R. Fry, J, G\'{o}mez, J.A. Jaramillo, M. Lovo, {\em
$C^1$-fine approximation of functions on Banach spaces with
unconditional basis}, Quarterly J. Math. 56 (2005), no. 1, 13--20.

\bibitem{DGZ2}
R. Deville, G. Godefroy, and V. Zizler, {\em A smooth variational
principle with applications to Hamilton-Jacobi equations in
infinite dimensions}, J. Funct. Anal. 111 (1993), no. 1, 197--212.

\bibitem{DGZ} R. Deville, G. Godefroy, and V. Zizler, {\em Smoothness
and renormings in Banach spaces}, vol. 64, Pitman Monographs and
Surveys in Pure and Applied Mathematics, 1993.

\bibitem{GaJaRan}
M. I. Garrido, J.A. Jaramillo, Y. Rangel, {\em personal
communication.}

\bibitem{Greene5}
R. E. Greene, and H. Wu, {\em On the subharmonicity and
plurisubharmonicity of geodesically convex functions}, Indiana
Univ. Math. J. 22 (1972/73), 641--653.

\bibitem{Greene3}
R. E. Greene, and H. Wu, {\em $C\sp{\infty }$ convex functions and
manifolds of positive curvature}, Acta Math. 137 (1976), no. 3-4,
209--245.

\bibitem{Greene4}
R. E. Greene, and H. Wu, {\em $C\sp{\infty }$ approximations of
convex, subharmonic, and plurisubharmonic functions}, Ann. Sci.
\'{E}cole Norm. Sup. (4) 12 (1979), no. 1, 47--84.


\bibitem{Klingenberg}
W. Klingenberg, {\em Riemannian geometry}, de Gruyter Studies in
Mathematics, 1. Walter de Gruyter \& Co., Berlin-New York, 1982.

\bibitem{LasryLions}
J.-M. Lasry, and P.-L. Lions,  {\em A remark on regularization in
Hilbert spaces.} Israel J. Math. 55 (1986), no. 3, 257--266.

\bibitem{Moulis} N. Moulis, {\em Approximation de fonctions diff\'{e}%
rentiables sur certains espaces de Banach,} Ann. Inst. Fourier,
Grenoble 21 (1971), no. 4, 293-345.

\bibitem{Myers}
S. B. Myers, {\em Algebras of differentiable functions}, Proc.
Amer. Math. Soc. 5, (1954). 917--922.

\bibitem{Nakai}
M. Nakai, {\em Algebras of some differentiable functions on
Riemannian manifolds}, Japan. J. Math. 29 1959 60--67.

\end{thebibliography}
\end{document}